# Star and Semi-Star Transformations in Fullerene Graphs

Meysam Taheri-Dehkordi

*University of Applied Science and Technology (UAST), Tehran, IRAN*

**m.taheri@uast.ac.ir**

## Abstract

A perfect star packing in a graph G is a spanning subgraph of G whose every component is isomorphic to the star graph $K_{1,3}$. A perfect star packing of a fullerene graph $G$ is of type $P0$ if all the centers of stars lie on hexagons of $G$. Many fullerene graphs arise from smaller fullerene graphs by applying some transformations. In this paper, we introduce two transformations for fullerene graphs that have perfect star packing of type $P0$ and examine some characteristics of the graphs obtained from this transformation.

## 1 Introduction

Fullerene graphs are cubic, 3-connected, planar graph with only pentagonal and hexagonal faces. By Euler's formula, it follows that the number of pentagonal faces in every fullerene graph is always twelve. In a paper by Grˊunbaum and Motzkin [8] showed that fullerene graphs with n vertices exist for all even n ≥ 24 and for n = 20. Klien and Liu [11] show that there exist fullerene graphs on n vertices with isolated pentagons for n = 60 and for each even n ≥ 70. For a systematic introduction on fullerene graphs, we refer the reader to the monograph [7] and for some symmetry-related aspects to [10].

A set of edges of graph G is said matching such that no two edges of it have a vertex in common. A matching is said to be perfect if every vertex of G is incident with an edge from its. In chemistry, the perfect matching is called Kekulé structure. In other words, a perfect matching is a spanning subgraph whose all components are isomorphic to K2. Suppose that G and L be two graphs. A spanning subgraph of G whose all components are isomorphic to L is called perfect L −packing in G. If L is the star graph K1,3, is called a perfect star packing. In [4], investigated which fullerene graphs allowed perfect star packings, and a type of perfect star packing was defined, which is called P0. As noted, many fullerene graphs arise from smaller fullerene graphs by applying some transformations. The number of vertices in the resulting graph is usually a multiple of the number of vertices of starting graph. The best known such transformation is the leapfrog transformation than can be thought of as a truncation of the dual. The number of vertices is tripled by this operation.

Two other operations, known as chamfer and capra, resulting in fullerene graphs with four and seven times, respectively, more vertices than the original graph. To study more information on these transformations, we refer the reader to [2, 10] and [1, 3]. In this paper we look at fullerene graphs having a perfect star packing of type P0 and define two transformations for these fullerene graphs.

## 2 Star Transformation

Let $G$ be a fullerene graph. A perfect star packing in $G$ is of type $P0$ if no center of a star is on a pentagon of $G$. From results of [5], such packings cannot exist in small fullerenes and only fullerenes with isolated pentagons can have such perfect star packings.

Theorem [5]. *The smallest fullerene having a perfect star packing of type $P0$ is the icosahedral isomer of $C_{80}$.*

Using the definition of this packing we introduce a transformation for fullerene graphs. We call this transformation, the star transformation.

Let $G$ be a fullerene graph with a perfect star packing of type $P0$. Let $u$ be in the center of one star, as shown in Figure 1.

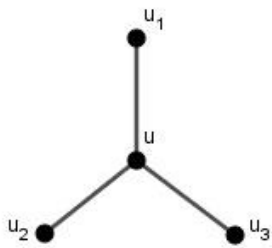

**Figure 1.** $u$ is in the center of one star

We connect each pair of vertices $u_1$, $u_2$, and $u_3$. We do this for all stars in the perfect star packing in $G$. (See Figure 2).

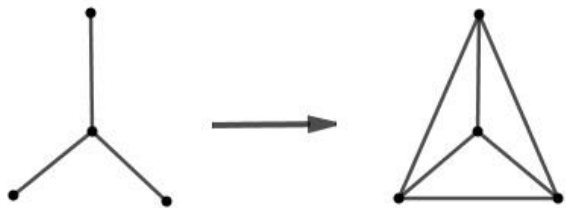

**Figure 2.** Connecting vertices adjacent to the central vertex of a star

By continuing this process, we will have $\frac{n}{4}$ subgraphs as shown in Figure 2 left.

We subdivide each new edge by one vertex and connect each of these added vertices to its corresponding vertex in the same hexagon. (Dashed lines in Figure 3.)

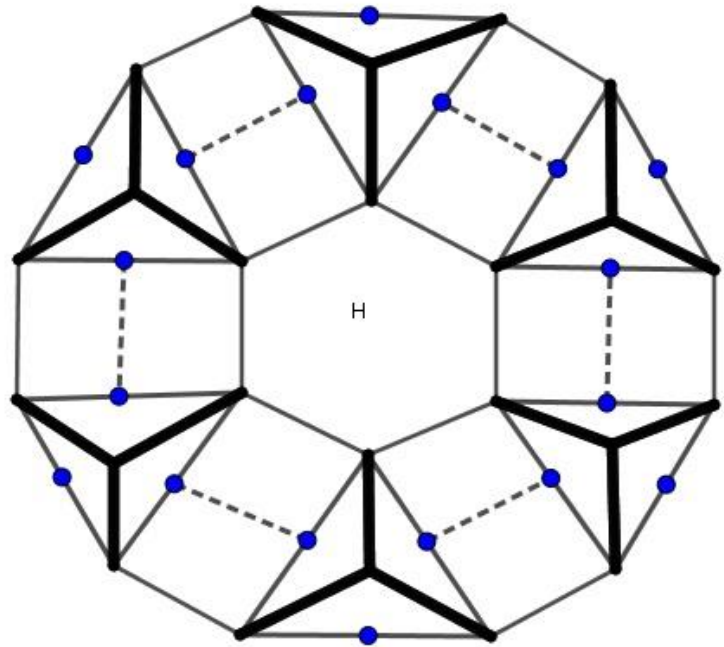

**Figure 3.** How the star transformation works

As we know, in a perfect star packing of type $P0$, the center of each star is shared by three hexagons. We consider hexagons that none of whose vertices is the central vertex of the star, (Like the hexagon in Figure 3). Into each of these hexagons and also in all pentagons, inscribe a polygon with the same number of sides. Connect each vertex of the original fullerene with three vertices of newly inscribe polygons. (See Figure 4)

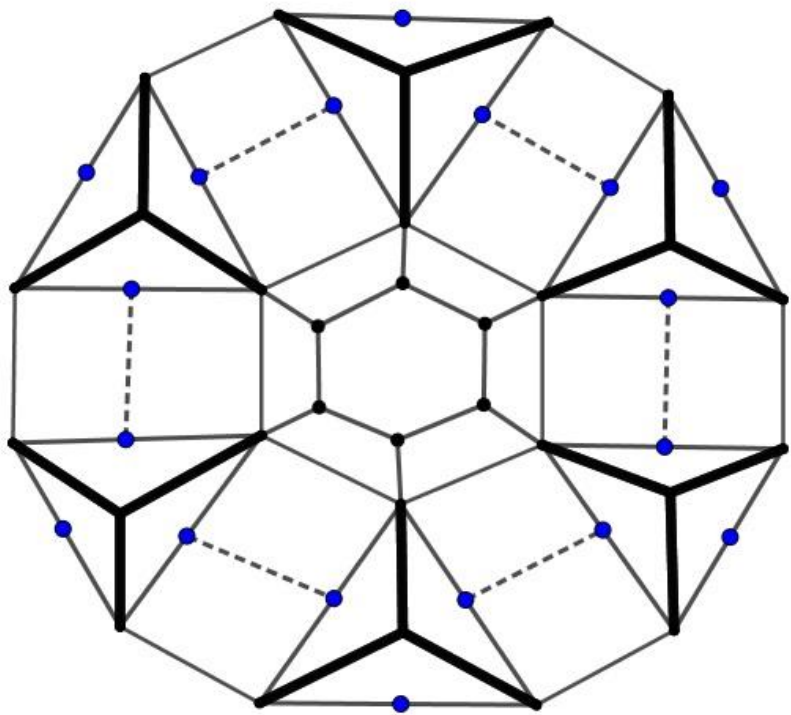

**Figure 4.** Connecting each vertex of the original fullerene with vertices of newly polygons

Finally, erase the center of stars in $G$. The faces of resulting graphs has only pentagons and hexagons. The number of pentagons in it is exactly twelve and the resulting graph is planar, cubic and 3-connected and hence it is a fullerene graph.

Proposition. *The number of vertices of the fullerene graph created by the star transformation is $\frac{9}{4}$ times the number of vertices in original fullerene.*

Proof. Let $F$ be a fullerene graph arising from the fullerene graph $G$ via the star transformation. Suppose $u$ is in the center of one star in $G$ and $u_1, u_2$ and $u_3$ be vertices connected to $u$. Via the Star transformation, Vertices $v_1, v_2, \ldots, v_6$ are created. So, the four vertices $u, u_1, u_2,$ and $u_4$ become nine vertices $u_1, u_2, u_3, v_1, v_2, v_3, v_4, v_5$ and $v_6$. (See Figure 5) Therefore, the number of vertices of $F$ is $\frac{9}{4}$ times the number of vertices of $G$.

Finally, because $G$ has a perfect star packing, the number of its vertices is a multiple of four, and therefore $\frac{9}{4} n(G)$ is always an integer.

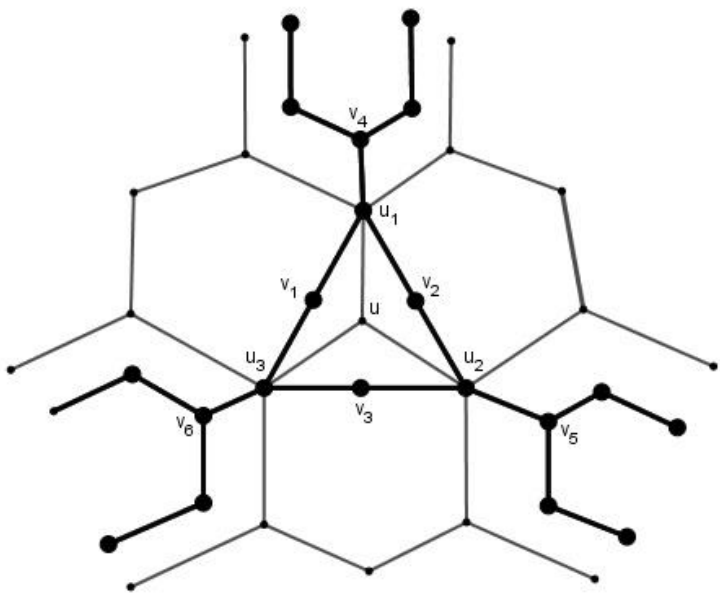

**Figure 5.** Added vertices in star transformation

The effect of star transformation on the graph with 80 vertices and the resulting graph with 180 vertices is shown in Figure 6.

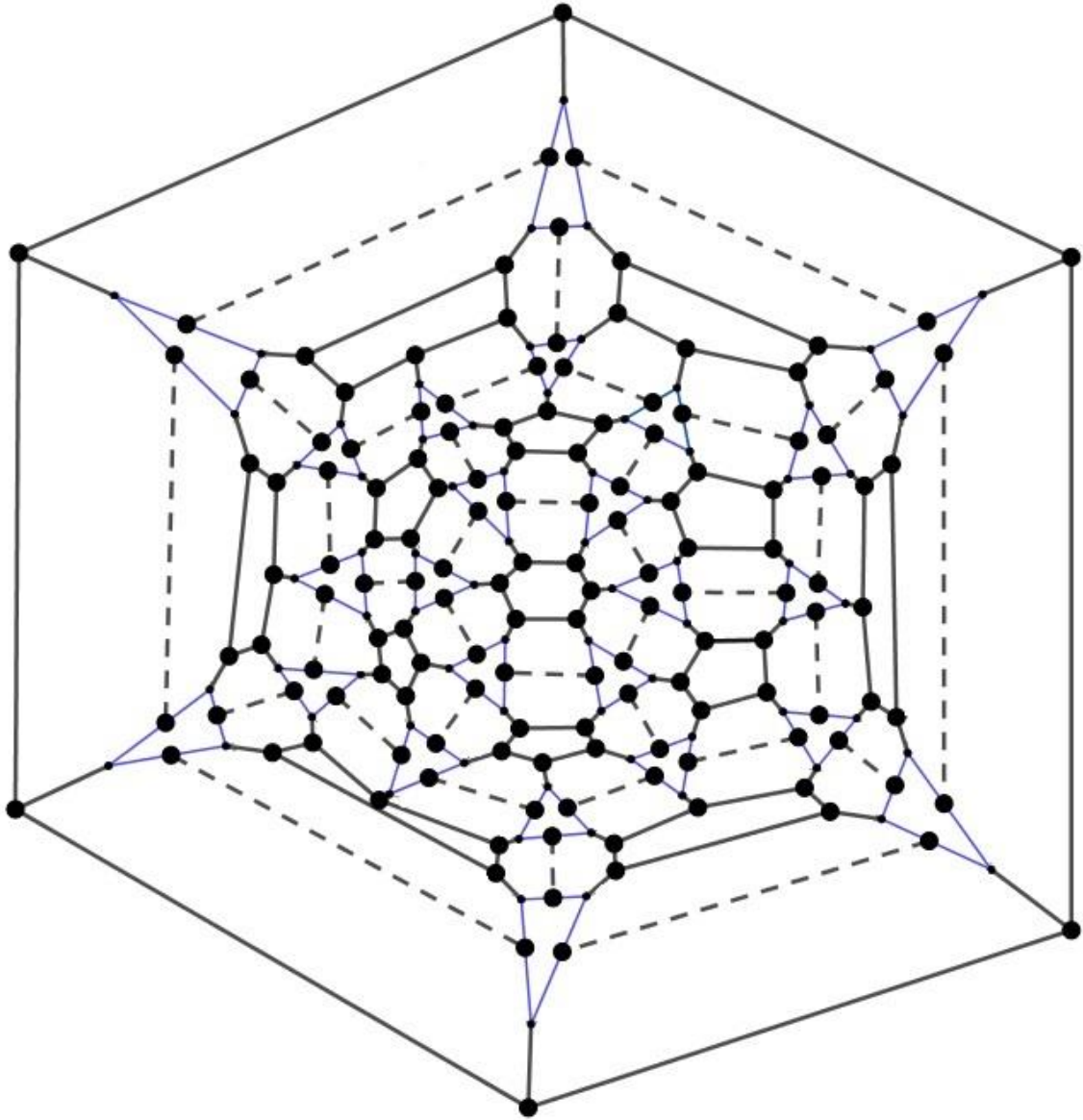

**Figure 6.** A graph of 180 vertices created by star transformation

From paper [5], we know that the smallest fullerene graph having a perfect star packing of type $P0$ has 80 vertices.

Table 1, shows the number of vertices of some fullerene graphs $F$ arising from a fullerene graph $G$ via the star transformation.

|     | $n(G)$ | $n(F)$, Star Transformation |
| --- | --- | --- |
| 1   | 80  | 180 |
| 2   | 96  | 216 |
| 3   | 104 | 234 |
| 4   | 112 | 252 |
| 5   | 120 | 270 |
| 6   | 128 | 288 |
| 7   | 136 | 306 |
| 8   | 144 | 324 |
| 9   | 152 | 342 |
| 10  | 160 | 360 |

**Table 1.** The number of vertices via star transformation

From paper [13], we know that only fullerene graphs with $8n$ vertices have a perfect star packing. Therefore, in the above table, only fullerene graphs of rows 2, 6 and 10 can have perfect star packing. Therefore, as an open problem, it is possible to check which of the graphs with $8n$ vertices arising from star transformation have star packing.

Theorem. *Let $G$ be a fullerene graph on n vertices, $n \geq 80$, and $F$ be a fullerene graph arising from $G$ via the star transformation. Then $F$ has a perfect $\{C_5, C_6\}$ – packing.*

Proof. According to how the star transformation works, all vertices of $F$ can be covered by a collection of disjoint cycles of sizes 5 and 6.

A perfect pseudo matching in a graph $G$ is a spanning subgraph of $G$ whose components are isomorphic to $K_{1,3}$ or to $K_2$. For more on the problem of pseudo matching and the size of a perfect pseudo matching in fullerene graphs, see [5,6].

Theorem. *Let $G$ be a fullerene graph on $n$ vertices and $F$ be a fullerene graph arising from $G$ via the star transformation. Then $F$ has a perfect pseudo matching with two components isomorphic to star graph $K_{1,3}$.*

Proof. Consider a hexagon adjacent to five pentagons, in $F$. We cover the vertices of this hexagon as shown in Figure 7. Around this hexagon, there are six other hexagons, which are covered in the form of Figure 7. Other vertices can be covered by components isomorphic to $K_2$.

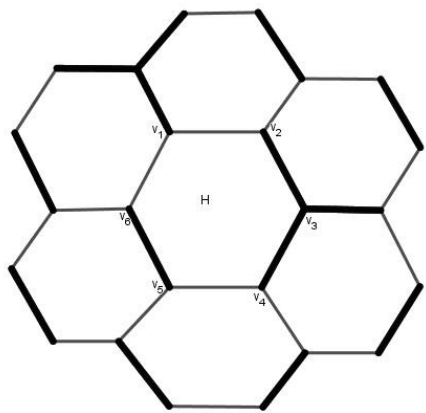

**Figure 7.** Packing of hexagon $H$

An example of this type of packaging is shown in Figure 8.

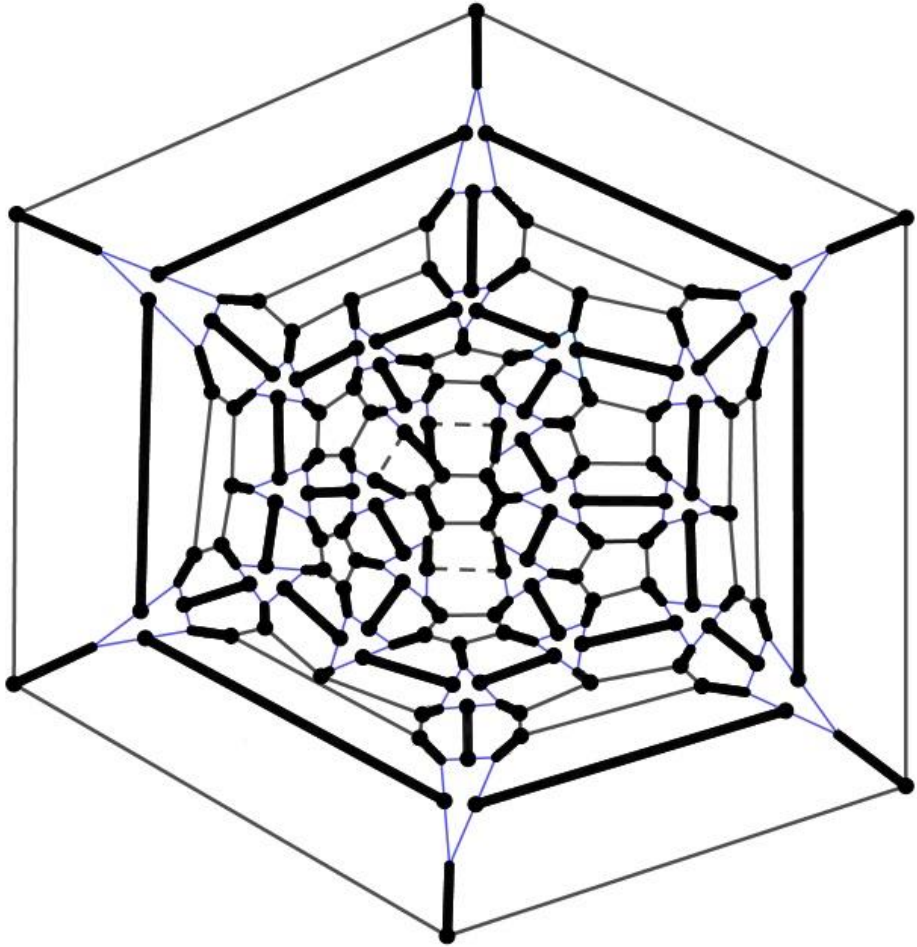

**Figure 8.** Packaging with two stars (Thick lines)

Theorem [3]. *Let G be a planar, 3-regular and 3-connected graph whose faces have a maximum length of six. Then G is Hamiltonian.*

From the previous theorem, we know that all fullerene graphs are Hamiltonian.

Theorem. *Let G be a fullerene graph with a perfect star packing of type P0 and F arises from G via the star transformation. Then there is a perfect packing of $P_3 = K_{1,2}$ and $P_9$ in G.*

Proof. Considering that the number of vertices of $F$ is a multiple of 9 and according to the Hamiltonian of fullerene graphs, the Hamiltonian path in $F$ can be considered as follows.

$$v_1 v_2 v_3 v_4 v_5 v_6 v_7 v_8 v_9 v_{10} \cdots v_{9n-2} v_{9n-1} v_{9n}$$

According to the number of vertices, we divide the above path into paths of length 8. Therefore, a packing of $P_9$ is created. By a similar argument, the above Hamiltonian path can be divided into paths that are isomorphic to $P_3$. □

## 3 Semi-Star Transformation

Let $G$ be a fullerene graph on $n$ vertices ($n \geq 80$) with a perfect star packing of type $P0$ arising from a smaller fullerene via the star transformation. Then $G$ has an even number of the star graph $K_{1,3}$. All the vertices of graph $G$ can be covered by stars by packing as shown in Figure 9.

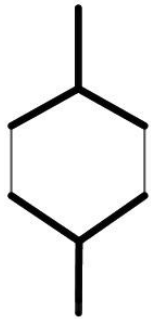

**Figure 9.** Packing a hexagon by two stars

Subdivide each edge of the stars by one vertex. As shown in Figure 10, we connect two vertices in all hexagons.

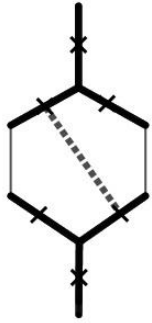

**Figure 10.** Connecting two vertices in a hexagon (Dash line)

By doing this process, we arrive at a fullerene graph $F$ which has a twice the number of hexagons.

According to the way $F$ is created, the number of vertices of $F$ is equal to

$$n + \frac{3n}{4} = \frac{7n}{4}$$

Let me call this transformation the **semi-star transformation**.

The effect of semi-star transformation on a fullerene graph on 80 vertices with a perfect star packing of type $P0$ is shown in Figure 11.

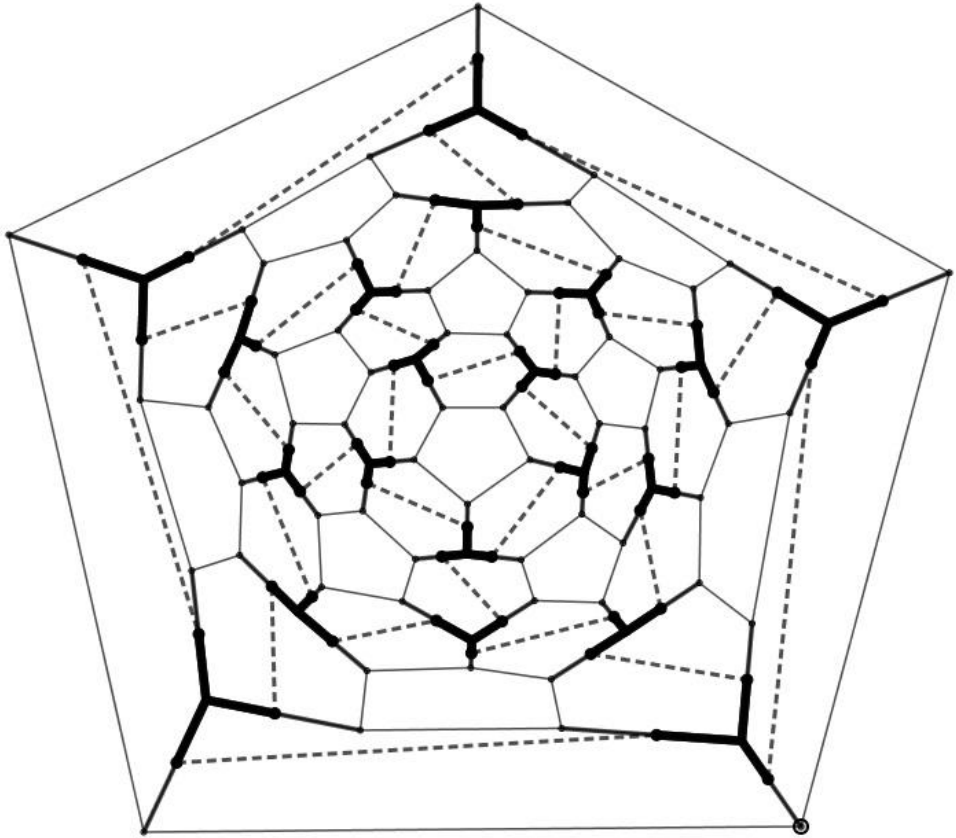

**Figure 11.** Semi-star transformation in $C_{80}(I_h)$

**Theorem.** *Let $G$ be a fullerene graph on $n$ vertices with a perfect star packing of type $P0$ arising from a smaller fullerene via the chamfer transformation. If $F$ arises from $G$ via semi-star transformation, then $F$ has a perfect packing of $S(K_{1,3})$.*

**Proof.** Consider an arbitrary star in the perfect packing star of $G$. According to the act of the star transformation, if we create a vertex on each edge of the stars, a graph will be created as shown in Figure 12. which is a subdivided star that covers all vertices of $F$.

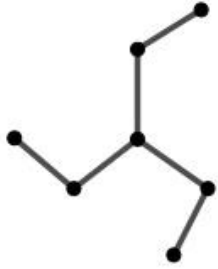

**Figure 12.** Subdivided star $S(K_{1,3})$.

One of the subjects of interest is whether the graphs created by star transformation have star packing or not.